\input harvmac
\null
\vskip 2truecm

\font\tit=cmr12
\centerline{\tit KAPPA-CONTRACTION FROM $SU_q(2)$ TO
$E_{\kappa}(2)$\footnote{$^*$}{{Submitted to the Proceedings of the 4-th
Colloqium {\it Quantum Groups and Integrable Systems}, Prague June
1995}}
}
\vskip 1truecm
\centerline{J. Sobczyk}
\medskip
\centerline{Institute for Theoretical Physics, Wroc\l aw University, pl.
Maxa Borna 9}
\centerline{50-204 Wroc\l aw, Poland; e-mail:
jsobczyk@proton.ift.uni.wroc.pl}
\vskip 2truecm
\centerline{ABSTRACT}
We present contraction prescription of the quantum groups: from
$SU_q(2)$ to $E_{\kappa}(2)$. Our strategy is different then one chosen
in ref. [P. Zaugg, J. Phys. A {\bf 28} (1995) 2589]. We provide
explicite
prescription for contraction of $a, b, c$ and $d$ generators of
$SL_q(2)$ and arrive at $^*$ Hopf algebra $E_{\kappa}(2)$.
\vskip 4truecm

The study of deformations of the $E(2)$ - two-dimensional
euclidean
group (or its dual $U(e(2))$) reveals many interesting features of the
theory of quantum groups.
It turns that they admits different
deformations according to different Lie bialgebra structures on $e(2)$.
One of the most
fruitful (and historically the first) approaches to contruct quantum
deformations of $U(e(2))$ is that of contraction from $U_q(sl(2))$ or
$SL_q(2)$
\ref\cgsta{E. Celeghini, R. Giachetti, R. Sorace and M. Tarlini, J.
Math. Phys. {\bf 31} (1990) 2548; {\it Contractions of quantum groups}
in {\it Quantum Groups}, Lecture Notes in Mathematics 1510, 221,
Springer Verlag 1992} -
\nref\wor{S. Woronowicz, Lett. Math. Phys. {\bf 23} (1991) 251}
\ref\swz{P. Schupp, P. Watts and B. Zumino, Lett. Math. Phys. {\bf 24}
(1992) 637}.
In this
way two different deformations of $U(e(2))$ were discovered \cgsta . It
is
also known that with the help of contraction it is possible
to obtain one of the two dual quantum deformations of
$E(2)$ \swz . With the second quantization of $E(2)$, as far as the
contraction is concerned, the situation is less clear. This contraction
involves the renormalization of the deformation parameter and is a $D=2$
analogue of the $\kappa$-deformation of the $D=4$ Poincar\'e algebra.
One possible approach to the problem presented in ref.
\ref\zau{P. Zaugg, J. Math. Phys. {\bf 36} (1995) 1547; J. Phys. A {\bf
28} (1995) 2589}
takes
advantage of the
quantum plane and is essentially a contraction of $U_q(O(3))$ quantum
group. In this
contribution we would like to follow the more straithforward route and
present a contraction scheme in terms of $a, b, c$ and $d$ generators of
the quantum $SL_q(2)$ together with the $*$ operation making it
$SU_q(2)$. They are subject to the standard relations
\eqn\stan{T = \pmatrix{a&b\cr c&d}\qquad RT_1T_2 = T_2T_1R}
with
\eqn\rmat{R =\pmatrix{q&0&0&0\cr 0&1&0&0\cr 0&q-q^{-1}&1&0\cr 0&0&0&q}}
\eqn\stcop{\Delta (T_{ij}) = T_{ik}\otimes T_{kj} }
\eqn\copstar{S(T)=T^{-1},\qquad T_{jk}^* = S(T_{kj})}
We will show that the following prescription defines a required
contraction $SU_q(2)$ $\rightarrow$ $E_{\kappa}(2)$ in the limit as the
contraction parameter $\rho\rightarrow\infty$
\eqn\aa{a = K + {L\over\rho} + ...\qquad b = M + {iN\over
\rho}
+ ...}
\eqn\ab{c = M - {iN\over \rho} + ...\qquad q = e^{{1\over
\kappa\rho}}}

These formulas were suggested (but in rather indirect way) by the
duality relations between $Sl_q(2)$ and $U_q(sl(2))$ and the
contraction prescription $U_q(sl(2))$ $\rightarrow$ $U_{\kappa}(e(2))$.
We did not write down the contraction behaviour of $d$
as it is dictated by the determinant relation
\eqn\ac{ad - q bc = 1}
One can show:
\eqn\ad{ d=K^{-1}(1+M^2) + {1\over\rho}\left( {1\over
\kappa}K^{-1}M^2
- K^{-1}LK^{-1}(1+M^2)\right) + ...}
In order to find this expression one has to assume that operator $a$
(and then also $K$) is invertible. This is a typical situation in
contraction of quantum
groups. Some discussion related to this point can be found in \wor ,
\ref\soel{D. Ellinas and J. Sobczyk, J. Math. Phys. {\bf 36} (1995)
1404}. We
start now to analyze relations imposed on $K, L, M$ and $N$ by \stan\
and
\stcop . In the zero order in ${1\over\rho}$ \stcop\ give rise to
\eqn\ca{K^2 = 1 + M^2}
Next terms linear in ${1\over\rho}$ impose the condition:
\eqn\cb{[L, K] = {1\over \kappa} M^2}
Taking all that into account one can write coproducts for our set of
variables
\eqn\cc{\Delta (K) = K\otimes K + M\otimes M}
\eqn\cd{\Delta (M) = K\otimes M + M\otimes K}
\eqn\ce{\Delta (L) = K\otimes L + L\otimes K + iN\otimes M
- iM\otimes N}
\eqn\cf{\Delta (N) = K\otimes N - i L\otimes M + i M\otimes L + N\otimes
K}
The star operation inherited from $SL_q(2)$ acts on them as follows:
\eqn\sta{K^*=K,\qquad L^*=-L,\qquad M^*=-M,\qquad N^*=-N
-{iM\over\kappa}}

It is useful to introduce linear combinations of $K, L, M$ and $N$,
namely $K\pm M$ and $L \mp {M\over 2\kappa} \pm iN$.
Their coproducts are:
\eqn\da{\Delta (K\pm M) = (K\pm M)\otimes (K\pm M)}
\eqn\db{\Delta (L \mp {M\over 2\kappa} \pm iN) = (L \mp
{M\over 2\kappa} \pm iN)\otimes (K\pm M)
+ (K\mp M)\otimes (L \mp {M\over 2\kappa}\pm iN) }
Now we investigate what follows from the contraction of quadratic
relations given in \stan .
$ab = q ba$, again up to terms ${1\over\rho}$, implies
\eqn\ea{[K, M] = 0}
\eqn\eb{[K, iN] + [L, M] = {1\over\kappa} MK}
Next, $ac = qca$ implies
\eqn\ed{[L, M] - [K, iN] = {1\over\kappa} MK}
Altogether we obtain
\eqn\ee{[K, N] = 0}
\eqn\ef{[L, M] = {1\over\kappa} MK}
Finally, $bc = cb$ is transformed into
\eqn\eg{[M, N] = 0}

In the variables $K\pm M$ and $L \mp {M\over 2\kappa}\pm
iN$ we first observe that
\eqn\fa{(K + M)(K - M) = (K - M)(K + M
) = 1}
Then we calculate
\eqn\fb{[L\mp{M\over 2\kappa}\pm iN, K + M] = {1\over
2\kappa} \left( (K+M)^2 -1\right)}
\eqn\fbb{[L\mp{M\over 2\kappa}\pm iN, K - M] = {1\over
2\kappa} \left( (K-M)^2 -1\right)}

On the other hand we are yet unable to calculate the commutator
\eqn\fc{[L - {M\over 2\kappa} + iN, L + {M\over 2\kappa} -
iN]}
as we do not know what is $[L, N]$. From the form of \aa -\ab\ it is
clear that we really need the analysis of terms of
order
${1\over\rho^2}$. We shall come back to this point in the discussion.
Let us introduce still another set of variables
\eqn\ga{\eta = (L - {M\over 2\kappa} + iN)(K - M)}
\eqn\gb{\bar\eta = - (K + M)(L + {M\over 2\kappa} - iN)}
They are defined in such a way that $\eta^* = \bar\eta$. On the other
hand $(K + M)^* = K - M$.  Let us call $(K + M)^2 = e^{i\alpha}$. We
derive
\eqn\gc{\Delta (\eta ) = \eta\otimes 1 + e^{-i\alpha}
\otimes\eta}
\eqn\gd{\Delta (\bar\eta ) = \bar\eta\otimes 1 + e^{i\alpha}
\otimes\bar\eta}
\eqn\ge{\Delta (e^{i\alpha}) = e^{i\alpha}\otimes e^{i\alpha}}
\eqn\gf{ [\eta, e^{i\alpha}] = {1\over\kappa} (e^{i\alpha} -1)}
\eqn\gf{ [\bar\eta, e^{i\alpha}] = {1\over\kappa}
(e^{i\alpha} - e^{2i\alpha})}
As we have said before, for completeness we still need the commutator
$[\eta, \bar\eta ]$. If we look at the coproducts of $\eta$ and
$\bar\eta$ one discovers that the consistent choice is
\eqn\gg{ [\eta, \bar\eta ] = {1\over\kappa} (\bar\eta + \eta )}
It is interesting that a similar problem of a necessity of completing
the structure of the quantum group appeared before in the paper
\ref\bcgst{A. Ballesteros, E. Celeghini, R. Giachetti, R. Sorace and M.
Tarlini, J. Phys. A {\bf 26} (1993) 7495}
where the authors tried to determine the quantum releations of $E_q(2)$
by means of the numerical $R$ matrix satisfying Yang Baxter equation.
In our case the remaining relation can certainly be obtained by analysis
of the
terms quadratic in the contraction parameter. This analysis is however
very tedious and not particulary iluminating. Let us still mention that
expressions for the antipode and counit can be also obtained along the
lines explained above.
In this way we reproduced the relations defining the quantum group
$E_{\kappa}(2)$. Expressions \gc - \gg\ can be compared to those
obtained before with a help of other techniques e.g in refs. \bcgst\
\ref\mas{P. Ma\'slanka, J. Math. Phys. {\bf 35} (1994) 1976}.
\listrefs
\end